\documentclass{amsart}

\usepackage{amsmath,enumerate}
\usepackage{amsthm}
\usepackage{amssymb}
\usepackage{amsfonts}
\usepackage{hyperref}

\newtheorem{theorem}{Theorem}[section]
\newtheorem*{theorem*}{Theorem}
\newtheorem{proposition}[theorem]{Proposition}
\newtheorem{lemma}[theorem]{Lemma}
\newtheorem{corollary}[theorem]{Corollary}

\newcommand{\rdots}{\mathinner{%
  \mkern1mu\raise1pt\hbox{.}%
  \mkern2mu\raise4pt\hbox{.}%
  \mkern2mu\raise7pt\vbox{\kern7pt\hbox{.}}\mkern1mu}}

\newcommand{\tmop}[1]{\operatorname{#1}}

\begin{document}

\title{A Casselman-Shalika formula for the Shalika model of $\tmop{GL}_n$.}
\author{Yiannis Sakellaridis} 

\address{Department of Mathematics \\
Stanford University \\
Stanford, CA 94305-2125 \\
USA}

\email{yiannis@math.stanford.edu}

\begin{abstract}
The Casselman-Shalika method is a way to compute explicit formulas for periods of irreducible unramified representations of $p$-adic groups that are associated to unique models (i.e. multiplicity-free induced representations). We apply this method to the case of the Shalika model of $\tmop{GL}_n$, which is known to distinguish lifts from odd orthogonal groups. In the course of our proof, we further develop a variant of the method, that was introduced by Y.Hironaka, and in effect reduce many such problems to straightforward
  calculations on the group.
\end{abstract}

\subjclass[2000]{Primary 22E50; Secondary 11F70,11F85}

\maketitle


\section{Introduction}

Let $G=\tmop{GL}_{2n}(F)$, where $F$ is a non-archimedean
local field. Let $S$ denote the ``Shalika subgroup'', consisting
of matrices of the form 
\begin{equation}\left( \begin{array}{cc}
 g \\ & g \end{array}\right) \cdot
\left(
\begin{array}{cc} I & X \\ & I \end{array}\right)
\end{equation}
with $g\in \tmop{GL}_n(F)$ and $M\in \tmop{Mat}_n(F)$, and let
$\psi$ be an additive character  of $F$
whose conductor is the ring of integers $\mathfrak o$. The expression
\begin{equation}\psi(\tmop{tr}(X))\end{equation}
defines a character $\Psi$ on the Shalika subgroup. By a Shalika model 
for a smooth representation 
$\pi$ of $G$ we mean a $G$-equivariant
 morphism from the space of $\pi$ into the representation induced from 
the character $\Psi$ on $S$, which can be realized as the space of 
functions
 $f$ on $G$ such that:
\begin{equation}
f(sh)=\Psi(s)f(h)
\end{equation}
for all $s\in S,h\in G$.

It is known from the work of Jacquet and Rallis (\cite{JR}) that
every irreducible admissible representation $\pi$ of $\tmop{GL}_n$
over a non-archimedean local field\footnote{Jacquet and Rallis require
that the field is of characteristic zero, although it is not clear to me
whether they need this requirement. In any case, one can prove uniqueness specifically for the unramified principal series, and for characters in generic position, by a double coset decomposition $P\backslash G/S$, where $P$ is
a Borel subgroup (an extension of our arguments in section \ref{secvanishing}). Uniqueness for characters in generic position is enough for our purposes (it gives the correct formula for a spherical function in the Shalika model, even if it is not everywhere unique), therefore
we will not impose this condition on the field.} possesses at most one (up to
scaling by a constant factor) Shalika model. (Uniqueness at the
archimedean places is proven by Ash and Ginzburg in \cite{AG}.)
Furthermore, it is proven in the same paper that the existence of
those implies that $\pi$ is self-contragredient.

In this article we consider irreducible spherical representations (i.e. possessing a
vector invariant under the maximal compact subgroup), which we 
realize as (irreducible factors of) unramified principal series - i.e. induced
from an unramified character of the Borel subgroup of upper triangular
matrices.
Although I am not aware of any proof in print, it is probably
known for a long time that a principal series of $\tmop{GL}_n$ over
a local field admits a Shalika functional if and only if\footnote{With 
some additional complications for the degenerate - i.e. reducible - case,
which we explain later.} the induction
data is of the
form $(\chi_1, \chi_2, \dots, \chi_n, \chi_n^{-1},\dots,
\chi_2^{-1},\chi_1^{-1})$ (or a permutation of those).
This is the
image of lifts from $\tmop{SO}_{2n+1}$\footnote{Notice that although 
supercuspidal
lifts as well as
global automorphic lifts from $\tmop{SO}_{2n}$ and
$\tmop{SO}_{2n+1}$ are conjecturally disjoint, the spherical lifts are
not disjoint and therefore a principal series as above could also
be a lift from an even orthogonal group.}. This existence theorem
also follows from our arguments here (cf. sections
\ref{secvanishing} and \ref{secanalytic}). 

Our goal in this paper is to compute an explicit formula for the
spherical vector in the Shalika model of $I(\chi)$. For $n=1$ the
Shalika model coincides with the Whittaker model and the result is
well-known. For $n=2$ we have an isomorphism
$\tmop{PGL}_4=\tmop{PGSO}_6$ and the Shalika model coincides with
the ``Whittaker-orthogonal'' model for $\tmop{PGSO}_6$ considered
by Bump, Friedberg and Ginzburg in \cite{BFG}, therefore our
result is also known in this case. The result for $n\ge 3$ is new.
However, in the degenerate case where the spherical decomposition
factor of the representation is induced from a character of the
maximal parabolic with Levi factor $GL_n \times GL_n$, the spherical vector in the Shalika model for $F=\mathbb Q_p$ has a very interesting form,
which has been computed by F.Sato (\cite{Sato});  its value on elements of
the form
\begin{equation*}
\left(\begin{array}{cccc|c}
p^{\lambda_1} &&&&\\
& p^{\lambda_2} &&&\\
&& \ddots &&\\
&&&p^{\lambda_n} &\\
\hline &&&& I
\end{array}\right)
\end{equation*}
is equal to the product of an elementary factor with a 
certain generating function for subgroups of the finite abelian
p-group $\mathbb{Z}/p^{\lambda_1} \times \mathbb{Z}/p^{\lambda_2} \times
\cdots \times \mathbb{Z}/p^{\lambda_n}$. In the case $n=2$ the latter
was also proven by Bump and Beineke (\cite{BB}).

Our formula for the spherical vector appears in section
\ref{secprelim}, after some necessary definitions have been given.

Shalika models first appeared in the work of Jacquet and Shalika
(\cite{JS}); they constructed a Rankin-Selberg integral for
the exterior-square L-function:
\begin{equation}
L(s,\pi,\wedge^2\rho)
\end{equation}
where $\pi$ is an irreducible unitary cuspidal
representation of $\tmop{GL}_{2n}(\mathbb{A}_F)$ ($F$ now
a global field) and  $\rho$ is the standard representation of
$\tmop{GL}_{2n}({\mathbb C})$.
They showed that
this L-function has at most a simple pole at $s=1$, whose residue
is given by:
\begin{equation}
\int_{[\tmop{PGL}_{n}]}\int_{[\tmop{Mat}_n]}  \phi\left(\left(
\begin{array}{cc} g \\ & g \end{array}\right) \cdot \left(
\begin{array}{cc} I & X \\ & I \end{array}\right)h\right)
\psi(\tmop{tr}(X)) dX dg.
\end{equation}
(Here brackets denote the quotient of adelic points modulo $F$-points,
$\phi$ is a vector in the space
of $\pi$, $\psi$ a character of the adeles and we assume a trivial central 
character.)
This integral clearly represents a global Shalika functional.
Therefore,
the exterior-square L-function has a pole at $s=1$ if and only if
$\pi$ admits a global (nonzero) Shalika model.

Subsequently, Ginzburg, Rallis and Soudry proved in \cite{GRS}
that an automorphic representation induced from cuspidal data
$\pi_1,\pi_2,\dots,\pi_r$ is a (weak) lift from $\tmop{SO}_{2n+1}$
to $\tmop{GL}_{2n}$ (corresponding to an inclusion of the
corresponding -connected- L-groups: $\tmop{Sp}_{2n}(\mathbb
C)\to\tmop{GL}_{2n}(\mathbb C)$) if and only if the
exterior-square L-functions of all induction data have a pole at
$s=1$; in particular, for cuspidal representations the existence
of global Shalika models characterizes (weak) lifts from
$\tmop{SO}_{2n+1}$.

In a different direction, Bump and Friedberg constructed in
\cite{BF} a Rankin-Selberg integral for the product of two
L-functions:
\begin{equation}
L(s_1,\pi)L(s_2,\pi,\wedge^2 \rho)
\end{equation}
(with $\pi,\rho$ as above). The residue of
that integral at $s_2=1$ is given by:
\begin{equation} \label{fj}
\int_{[(\tmop{GL}_{n}/Z)\times(\tmop{GL}_{n}/Z)]}
\phi\left(\left(\begin{array}{cc} g_1\\ &
g_2\end{array}\right)\right)
\left|\frac{\tmop{det}(g_1)}{\tmop{det}(g_2)}\right|^{s_1-1/2}
dg_1 dg_2.
\end{equation}

The latter was examined by Friedberg and Jacquet in \cite{FJ}, who
showed that it is  nonzero for some $\phi$ in the space of $\pi$
if and only if $L(s,\pi,\wedge^2\rho)$ has a pole at $s=1$. In
that case, the latter integral unfolds to an integral involving
the global Shalika model for $\pi$, and for some $\phi$ it
represents the L-function $L(s_1,\pi)$.

The formula that we give could be useful in studying
Rankin-Selberg integrals which unfold to the Shalika model,
although we must note that the Rankin-Selberg integrals discussed
above have already been examined without use of such a formula. In addition, Shalika models
appear as Fourier coefficients of Eisenstein series on
$\tmop{GL}_{2n}$. (This is the motivation behind the
aforementioned work of Beineke and Bump (\cite{BB}) and of Sato (\cite{Sato})
on the degenerate case.) One could also use the formula to do explicit harmonic analysis on the space $\tmop{Ind}_S^G(\Psi)$. Finally, the methods that we use reduce the computation of such a formula to a straightforward calculation
on the group, and therefore their scope extends beyond the Shalika model.

For the proof, we follow the method initiated in Casselman
(\cite{C}) and Casselman and Shalika (\cite{CS}). The basic
ingredient there is to express the spherical vector in the
so-called ``Casselman basis'' of $I(\chi)^B$, i.e. invariants of
the standard Iwahori subgroup $B$. This basis is indexed by
elements of the Weyl group, and the final formulas are expressed
as sums over the Weyl group. In \cite{C} and \cite{CS} this method
is used to prove MacDonald's formula for the spherical functions
and an explicit formula (first proven by Shintani for
$\tmop{GL}_n$ and Kato for Chevalley groups) for the spherical
Whittaker function.

However, one runs into computational difficulties in other cases.
We explain them briefly in section \ref{secbasis}. Therefore,
attempts to apply this method to other models have typically
involved tedious calculations and arguments. An alternative
introduced by Y.Hironaka in \cite{H} provides a framework of much
wider applicability and smaller complexity. The basic new idea
here is to express in the Casselman basis, instead of the
spherical vector, a suitable projection of the functional which is
used to define the model. We also explain this in section
\ref{secbasis}.

Subsequently, this variant was used successfully by Mao and Rallis
in \cite{MR} and by Omer Offen in \cite{O} to examine certain
models of self-contragredient principal series, i.e.
representations like the ones that we, too, are considering here.
They found that some of the terms in the Casselman-Shalika formula
vanish in those cases - for instance in the work of Offen the
final formula is a sum over the Weyl group, not of
$\tmop{GL}_{2n}$, but of $\tmop{Sp}_{2n}$. This should be relevant
to the fact that the representation is a lift from
$\tmop{SO}_{2n+1}$, whose L-group is $\tmop{Sp}_{2n}(\mathbb C)$.

The arguments that we use are very close to the ones used in
\cite{MR} and \cite{O}. However, there is no explicit
consideration of a symmetric space, and at several points we have
used different arguments, which are easily applied 
in many different contexts, as we explain in section \ref{secgeneral}.

An outline of the paper is as follows:

In section \ref{secprelim} we introduce notation and state our
main theorem. In section \ref{secopen} we show that there exists
an open orbit of the Shalika subgroup in the flag variety
of $G$. This fact is closely related to the uniqueness of Shalika
models, but we don't expand more on this idea here. In section
\ref{secbasis} we review the Casselman-Shalika method and its
variant introduced by Hironaka. In section \ref{secvanishing} we
show (as in \cite{MR} and \cite{O}) that the only terms which
appear in the final formula are those which correspond to the Weyl
group of $\tmop{Sp}_{2n}$. We also show that the Shalika
functional, when applied to functions which are supported on the
open orbit, has a simple integral representation. In section
\ref{secdependence} we use this integral representation to perform
a simple calculation, and in section \ref{secanalytic} we show
that this integral representation is valid even for functions not
supported on the open orbit, provided that the inducing parameters
for $I(\chi)$ lie in some open subregion. We also use Bernstein's 
lemma to show that the
result should vary rationally with the induction data,
therefore allowing us to focus on the region of convergence only.
In section \ref{secfe} we complete the proof by computing
the effect of intertwining operators on the Shalika functional. Finally, section \ref{secgeneral} contains a discussion of how one might use our methods to compute (asymptotic, in general) values of the spherical vector in any unique model induced from a character of a closed algebraic subgroup of a split reductive group.

\textbf{Acknowledgements.} I thank Professor Daniel Bump for
suggesting the problem, for many useful discussions while I was
working on it, and for his help in preparing this article. This work was in part supported by NSF grant FRG DMS-0354662.


\section{Preliminaries and statement of the result}
\label{secprelim}

By $F$ we will denote a non-archimedean local field, by $\mathfrak o$ its
ring of integers, by $q$ the order of its residue field
and by $\varpi$ a uniformizing element. The group
$\tmop{GL}_{2n}(F)$ will be denoted by $G$ and its Borel subgroup of 
upper triangular matrices by $P$. We let $\chi=(\chi_1,\chi_2,\dots,\chi_{2n})$
denote an unramified character of $P$ and $I(\chi)$ the smooth unramified
principal series representation of $G$, obtained by (normalized)
induction from $\chi$. In other words, $I(\chi)$ consists of all smooth
(i.e. locally constant) functions on $G$ which satisfy:
\begin{equation}
f(pg)= \chi\delta^{\frac{1}{2}}(p) f(g)
\end{equation}
for every $p\in P, g\in G$, where 
\[ \chi\left(\left(\begin{array}{cccc}
a_1 & * & * & \cdots \\
& a_2 & * & \cdots \\
&& \ddots\\
&&& a_{2n}
\end{array}\right)\right) = \chi_1(a_1)\cdots
\chi_{2n}(a_{2n})\]
and $\delta=(|\cdot |^{2n-1},|\cdot |^{2n-3},\dots, |\cdot |^{-2n+1})$
is the modular character (the quotient of the
right and left Haar measures on $P$). 
The action of $G$ is by right translations: $R_g f (x) = f(xg)$. 

The Shalika subgroup $S$ consists of matrices of the form:
\begin{equation}s=\left( \begin{array}{cc}
 g \\ & g \end{array}\right) \cdot
\left(
\begin{array}{cc} I & X \\ & I \end{array}\right)
\end{equation}
and we use an additive character $\psi$ whose conductor is the ring of integers $\mathfrak o$ to define a character $\Psi(s)=\psi(\tmop{tr}X)$ on $S$. The space of the smooth induced representation $\tmop{Ind}_S^G(\Psi)$ consists of all locally constant functions $f$ that satisfy:
\begin{equation}
f(sg)=\Psi(s)f(g)
\end{equation}
for every $s\in S,g\in G$.

For any irreducible representation $\pi$, a $G$-equivariant morphism\footnote{We use the word ``model'' for such a morphism, although, strictly speaking, one should use the word ``model'' if the morphism is injective.}:
\begin{equation} \pi \to \tmop{Ind}_S^G (\Psi) \end{equation}
is equivalent to a ``Shalika functional'' $\Lambda$ on the space of $\pi$,
satisfying:
\begin{equation} \label{Shalika} \Lambda(\pi(s)v)=\Psi(s)\Lambda(v) \end{equation}
for every $s\in S, v$ in the space of $\pi$.
Indeed, the function $f_v(g)=\Lambda(\pi(g) v)$ will belong to  $\tmop{Ind}_S^G (\Psi)$
and, conversely, given such a morphism, the functor ``evaluation at $1\in G$''defines such a functional.

We will assume that $\chi$ is of the form:
\begin{equation}\label{chi}
\chi=(\chi_1,\chi_2,\dots,\chi_{n},\chi_n^{-1},\dots,
\chi_1^{-1})= (|\cdot|^{z_1},|\cdot|^{z_2}, \dots, |\cdot|^{z_n},
|\cdot|^{-z_n}, \dots, |\cdot|^{-z_1})
\end{equation}
or one of its translates through the action of the Weyl group
(i.e. permutations of the individual characters), for otherwise we
will see (section \ref{secvanishing}) that any Shalika functional on 
$I(\chi)$ is zero.

Let $K=\tmop{GL}_{2n}(\mathfrak o)$ be the standard maximal compact
subgroup of $G$ and $\phi_K$ the unique spherical (i.e. $K$-invariant)
vector in $I(\chi)$ normalized so that $\phi_K(1)=1$.
It is given by: $\phi_K(g)=
\chi\delta^{\frac{1}{2}}(p)$, where $g=p\cdot k$ is an Iwasawa
decomposition for $g$. Let $\Lambda$ denote a Shalika functional on 
$I(\chi)$, which is uniquely defined up to scaling (by the uniqueness
theorem of Jacquet and Rallis). Let
\begin{equation} \Omega(g)=\Lambda(R_g \phi_K) \end{equation}
denote the image of $\phi_K$ under the Shalika embedding defined by
$\Lambda$. For the moment, $\Omega$ is well defined only up to a scalar
factor, which may depend on $\chi$. 
\emph{Our goal is to compute an explicit formula for $\Omega(g)$.}

We notice first that it
suffices to compute it for a set of double coset representatives
in $S \backslash G /K$. By an easy argument (using the Iwasawa and
Cartan decompositions) we see that such a set of representatives
is given by the matrices:
\begin{equation}\label{glambda}
g_{\lambda}=\left(
\begin{array} {cc} \varpi^\lambda & \\ & I \end{array} \right )
\end{equation}
where $\varpi$ is a uniformizer for $F$, $\lambda=(\lambda_1,
\lambda_2, \dots,\lambda_n)$ with $\lambda_1\ge\lambda_2\ge\dots
\ge\lambda_n$ and $\varpi^{\lambda}$ denotes the matrix:
\[\left(
\begin{array}{cccc}
\varpi^{\lambda_1}\\
&\varpi^{\lambda_2}\\
&&\ddots\\
&&&\varpi^{\lambda_n}
\end{array}\right).\]

The statement of our formula involves the Weyl group of
$\tmop{Sp}_{2n}$ and some L-group terminology, therefore let us fix some notation:

Through the Satake isomorphism, the character $\chi$ is identified with
a semi-simple conjugacy class in the (connected) L-group $\tmop{GL}_{2n}
(\mathbb C)$, represented by the element $g_\chi = \tmop{diag}(\chi_1(\varpi), \chi_2(\varpi), \dots)$.
\emph{Notice that if the character is of the form (\ref{chi}) then $g_\chi$
can be considered as an element of $\tmop{Sp}_{2n}(\mathbb C)$}, i.e. the subgroup of $\tmop{GL}_{2n}(\mathbb C)$
stabilizing
the skew-symmetric bilinear form defined by
\[ J=\left(\begin{array}{cccccc} &&&&&1 \\
&&&&\cdots&\\
&&&1&&\\
&&-1&&&\\
&\cdots&&&&\\
-1&&&&& \end{array}\right).\] We denote by $W$ the Weyl group of $\tmop{GL}_{2n}$, by $\Gamma\simeq\mathbb{Z}/2\wr
S_n$ the Weyl group of $\tmop{Sp}_{2n}$, and by $\Phi_{GL},\Phi_{Sp}$ the corresponding root systems. By $\Phi_{Sp}^S$ we denote the set of short roots, by $\Phi_{Sp}^L$ the set of long roots, and a superscript $^+$ over $\Phi$
will denote positive roots (under the standard choice of those).

Our main theorem is:

\begin{theorem} \label{thm}
For a suitable normalization (given by (\ref{finalnorm})), the spherical Shalika function of
$I(\chi)$ is:
\begin{equation} \label{mainthm}
\Omega(g_\lambda)=\left\{\begin{array}{l}
\prod_{\alpha\in \Phi_{Sp}^+} (1-q^{-1}e^\alpha) \cdot \\
\,\,\,\,\,\,\,\,\,\,\,\, \cdot \delta^{\frac{1}{2}}(g_\lambda)\mathcal{A} \left(
e^{\rho+\lambda} \prod_{\alpha\in \Phi_{Sp}^{S+}}
(1-q^{-1}e^{-\alpha}) \right)(g_\chi) \mbox{                     , if }
\lambda_1\ge\dots\ge\lambda_n\ge 0 \\  0 \textrm{    ,otherwise}
\end{array}\right.
\end{equation}
where $\rho=\frac{1}{2}\sum_{\alpha\in \Phi_{Sp}^+}\alpha$ and
$\mathcal A$ denotes the ``alternator'':
\[ \mathcal{A}(\cdot)=\sum_{w\in \Gamma} (-1)^{l(w)} w(\cdot)\] ($l(w)$ the
length of $w$ in $\Gamma$).
\end{theorem}

This statement entails the claim that there exists a non-zero Shalika functional on $I(\chi)$ whenever the above expression is non-zero. This happens exactly when $\chi$ is regular (i.e. the matrix element $g_\chi$ is regular) and the spherical vector generates the principal series as a $G$-module (cf. section \ref{secfe}).

The reader who would like to avoid the L-group formalism should
take $e^\alpha(\chi)$ to mean $\chi(a_\alpha)$ where if $\alpha$
is the root $\tmop{diag}(t_1,t_2,\dots,t_{2n})\mapsto t_i
t_j^{-1}$ then $a_\alpha$ is the diagonal element with $\varpi$ on
the i-th line, $\varpi^{-1}$ on the j-th line and 1's otherwise;
and $e^\lambda(\chi)$ to mean $\chi(g_\lambda)$. Then, \emph{under
a different normalization} from above, the formula for the
spherical Shalika function (for $\lambda_1\ge\dots\ge\lambda_n\ge
0$) reads:
\begin{equation} \label{nolgroup}
\Omega(g_\lambda)= \sum_{w\in\Gamma} (-1)^{l(w)} \prod_{\alpha\in
\Phi_{Sp}^+,w\alpha<0} \chi(a_\alpha) \prod_{\alpha\in
\Phi_{Sp}^{S+},
w\alpha<0}\frac{1-q^{-1}\chi(a_{-\alpha})}{1-q^{-1}\chi(a_\alpha)}
{^w\chi}\delta^{\frac{1}{2}}(g_\lambda).
\end{equation}
The passage from one expression to the other is explained at the
end of section \ref{secfe}.

We shall also make use of the notation described below:

There exists a canonical surjection
$\mathcal {P}_{\chi}: C_c^{\infty}(G)\to I(\chi)$ given by:
\begin{equation}\label{Pchi} \mathcal{P}_{\chi} (f) (g)= \int_P
\chi^{-1}\delta^
{\frac{1}{2}} (p) f(pg) dp\end{equation} where the measure on $P$
will always be taken to be \emph{left} Haar measure. Under this
mapping, $\phi_K$ is just the image of the characteristic function
$1_K$ of $K$. Similarly, let $B$ denote the standard Iwahori
subgroup of $K$, consisting of matrices in $K$ which are upper
triangular modulo the prime $\mathfrak p$ of $F$, in other words,
whose entries below the diagonal belong to $\mathfrak p$. $K$ has
a Bruhat decomposition with respect to $B$: $K=\bigcup_{w\in W}
BwB$ (disjoint). We will denote $\phi_{BwB}=\mathcal P_{\chi}
(1_{BwB})$. When the character $\chi$ to which we are refering is
not obvious from the context, it will also appear as a subscript.

We denote by $A$ the maximal split torus of diagonal matrices, by $N$ the unipotent radical of the Borel subgroup, consisting of
upper triangular matrices with 1's on the diagonal, and by $N^-$ the opposite unipotent subgroup. For every root $\alpha\in \Phi$,
$N^{\alpha}$ will denote the image of standard embedding
corresponding to $\alpha$ of the additive algebraic group $\mathbb G_a \to
N$ (or $N^-$ if the root is negative), i.e. if $\alpha$ is the
root $t_i t_j^{-1}$ then $N^\alpha$ will consist of matrices with
1's on the diagonal and zeroes elsewhere, except for the ij-th
position. There exists a measure-preserving factorization:
\begin{equation}
N=\prod_{\alpha \in \Phi_{GL}^+} N^\alpha
\end{equation}
the product being taken in any order. The
image of $\mathbb G_a(\mathfrak o)$ under the above embedding will be denoted by $N^\alpha_0$ and
 the image of $\mathbb G_a(\mathfrak p^i)$ by $N^\alpha_i$. In general, a
subscript $_0$ will denote intersection with $K$, e.g. $P_0, N_0,
S_0$ etc. $N_1$ will denote the above product with $N^\alpha_1$'s
instead of $N^\alpha$, and $N_1^-$ will denote its transpose.  We
will also use $P^\alpha$ to denote that in the factorization
$P=A\cdot N$ the $N$ factor belongs to $N^\alpha$. Also, for a
simple root $\alpha$ we will use $N^{\widehat\alpha}$,
$P^{\widehat\alpha}$ to denote that the $N^\alpha$-factor is
missing (in other words, the $(i,j)$-th entry is 0, where now
$j=i\pm 1$). If $\alpha$ is a simple root, the simple reflection
corresponding to it will be denoted by $w_\alpha$. We will denote
the longest Weyl group element (both in $W$ and in $\Gamma$) by
$w_l$, and we shall identify elements in $W$ with permutation
matrices having only 1's and 0's as entries.

 The embedding
$\tmop{Sp}_{2n} \to \tmop{GL}_{2n}$ induces, dually, a
``collapse'' of the roots of $\tmop{GL}_{2n}$ to the roots of
$\tmop{Sp}_{2n}$ and identifies $\Gamma$ as a subgroup of the Weyl
group $W$ of $\tmop{GL}_{2n}$.  
The map:
\[ \Phi_{GL} \to \Phi_{Sp} \]
is one-to-one onto the set $\Phi_{Sp}^L$ of long roots of
$\tmop{Sp}_{2n}$ and two-to-one onto the set $\Phi_{Sp}^S$ of
short roots. If $\alpha$ and $\beta$ are two distinct roots in
$\Phi_{GL}$ that collapse to the same short root in $\Phi_{Sp}$, we
will write $\beta=\tilde{\alpha}$. For long roots, we adopt the
convention $\tilde\alpha=\alpha$.


\section{The open orbit} \label{secopen}

A functional $L$ on $I(\chi)$ corresponds to a distribution $D$ on
$G$ such that (by abuse of notation)
\begin{equation}\label{distr}
 D ( p g  ) = \chi^{- 1} \delta^{\frac{1}{2}} ( p ) D
( g
     )  \end{equation}
  for every $p \in P, g \in G$. The
  correspondence is given by:
 \begin{equation}  \label{corr} L (\mathcal{P}_{\chi} ( f  )) = D (f) =
     \int_{_{_{}} G} D ( x ) f ( x ) dx
     \end{equation}
  where $f \in C_c^{\infty} ( G )$, and we have used the usual integral
  notation for distributions.

This identifies the dual $I(\chi)^*$ of $I(\chi)$ with the space
of distributions satisfying (\ref{distr}). The \emph{smooth} dual
of $I(\chi)$ is the subspace $I(\chi^{-1})\subset I(\chi)^*$.

Let $\Delta$ be the distribution which corresponds to the Shalika
functional $\Lambda$. It satisfies the stronger relation:
\begin{equation}
\Delta ( p g s ) = \chi^{- 1} \delta^{\frac{1}{2}} ( p ) \Delta (
g
     ) \Psi^{- 1} ( s ) \end{equation}
  for every $p \in P, g\in G, s \in S$

Therefore $\Delta$ is fully determined by its ``values'' on a set
of representatives of
  double $P \backslash G / S$ cosets. We first prove:

  \begin{lemma} Let
\begin{equation}
\xi=\left(\begin{array}{cc}  & I \\ w_0&
  \end{array}\right)
\end{equation} where by $w_0$ we denote the $n\times n$ matrix \[\left(
\begin {array} {cccc} &&&1\\&&1&\\&\cdots&&\\1&&&
\end{array} \right).\]
Then the conjugate $H=\xi S \xi^{-1}$
  of $S$, has the property that $P\cdot H$ is Zariski open in $G$.

Equivalently, the open double coset $P\xi S$ is Zariski open.
\end{lemma}

\begin{proof}
Since
\[\xi\left(\begin{array}{cc}
g & X \\ & g \end{array}\right)\xi^{-1} = \left(\begin{array}{cc}
{g}& \\
{w_0 X}&{^{w_0}g} \end{array}\right)\] (the exponent on the left
denotes conjugation), the Lie algebra of $S$ consists of matrices
of the form
\begin{equation} \label{liealgebra}
\left(\begin{array}{cc}
{A}&\\
{B}&{^{w_0}A} \end{array}\right).\end{equation} It is then obvious
that it is complementary to the Lie algebra of the Borel subgroup,
hence the differential of the multiplication morphism $P\times H
\to G$ is surjective at the identity, therefore the image is
Zariski open.
\end{proof}

Since $G$ is irreducible as a variety, $P \xi S$ is the only open
double $(P\backslash, /S)$-coset.

Because of the above lemma, it will be more natural in most of the
proof to deal with the subgroup $H$ instead of $S$. Therefore, let
us see how things translate to this subgroup:

We consider the character $\Psi_H$ on $H$ defined by:
$\Psi_H(h)=\Psi(\xi^{-1}h \xi)$. There is a bijection
$\tmop{Ind}_S^G(\Psi) \to \tmop{Ind}_H^G(\Psi_H)$ given by $f
\mapsto f_H$ where
\begin{equation}\label{translate}
f_H(g)=f (\xi^{-1}g).
\end{equation}
Composing with the Shalika map we get a morphism
$I(\chi)\to \tmop{Ind}_H^G(\Psi_H)$. Let $\Lambda_H$ denote the
corresponding functional and $\Delta_H$ the corresponding
distribution, which satisfies:
\begin{equation}
\Delta_H ( p g h ) = \chi^{- 1} \delta^{\frac{1}{2}} ( p )
\Delta_H ( g
     ) \Psi_H^{- 1} ( h ) \end{equation}
  for every $p \in P, g\in G, h \in H$.

The spherical vector in the model induced from $H$ will be given
by:
\begin{equation}
\Omega_H(g)=\Omega(\xi^{-1} g)=\Omega(\xi^{-1}g\xi)
\end{equation}
(where we used the fact that $\xi \in K$), and in particular for
the representatives $g_{\lambda}$ as in (\ref{glambda}) we get:
\begin{equation} \label{omegaH}
\Omega(g_{\lambda})=\Omega_H (^\xi g_\lambda)= \Omega_H \left(
\left(
\begin{array} {cc} I & \\ & ^{w_0}\varpi^\lambda \end{array} \right )
\right) = \Omega_H \left( \left(
\begin{array} {cc} \varpi^{-\lambda} & \\ & I \end{array} \right )
\right).
\end{equation}

We caution the reader that, while we are looking at the
representatives $g_\lambda$ when refering to $\Omega$, we are
looking at $g_{-\lambda}$ when refering to $\Omega_H$.


\section{The Casselman basis} \label{secbasis}

In this section we summarize the method of Casselman and Shalika
and the variant of it that arises from the work of Hironaka. The
reader who is already familiar with this method and would like to skip this section should only keep
in mind that our goal in the rest of the paper will be to compute
the expression (\ref{goal}) which appears at the end of this
section.

The basic philosophy of the method is the following: Remember that
the behavior of the distribution $\Delta_H$ is determined modulo $P$
on the left and $H$ on the right. The expression
$\Omega_H(g_{-\lambda})$ involves the behavior of the distribution
$\Delta_H$ on the set $Kg_\lambda$, which intersects many
$P\backslash G /H$ double cosets, and as such is difficult to
handle. On the other hand we show that for functions $\phi\in
I(\chi)$ supported in $P\cdot H$, the Shalika functional has the
simple integral expression:
\begin{equation}\label{integral}
\Lambda_H(\phi)= \int_{H'} \phi(h) \Psi^{-1}(h) dh
\end{equation}
where $H'$ is some quotient of $H$.

Now, the computation is carried out by exploiting two facts:

1) It happens that for some $w\in W$, $BwBg_\lambda \subset
PH$
 for all $\lambda$. This allows us to compute the effect of the Shalika
functional on $R_{g_{-\lambda}}\phi_{BwB}$ by using (\ref{integral}).

2) The symmetries of $I(\chi)$ allow us to extend the computation
to all other Iwahori-invariants, i.e. elements of $I(\chi)^B$.

The second point certainly needs some clarification (and besides,
is only true in a very rough sense):

By ``symmetries'' we mean the fundamental fact that for $\chi$ in
general position (which means that the numbers $\pm z_i\pm
\frac{1}{2}$ are all distinct and hence $I(\chi)$ is irreducible),
$I(\chi)$ is isomorphic to $I(^w \chi)$ for every $w\in W$. This
is demonstrated by the intertwining operators $T_w: I(\chi)\to
I(^w \chi)$, which are $G$-equivariant maps (unique up to
scaling); for elements of $I(\chi)$ with support in
$\cup_{w'\nless w^{-1}}Pw'P$ they are given by the integral:
\begin{equation}
T_w \phi (1) = \int_{wNw^{-1}\cap N\backslash N} \phi(w^{-1}n) dn
= \int_{\prod_{\alpha>0, w^{-1}\alpha<0}N^\alpha} \phi(w^{-1}n) dn.
\end{equation}
(Remember that for us $w$ is represented by a permutation matrix.)

The connection between intertwining operators and
Iwahori-invariants arises from the fact that the operators are
``dual'' to $I(\chi)^B$ in a natural way: If we consider the
functionals on $I(\chi)$ defined by $\phi \mapsto T_w(\phi) (1)$,
($w\in W$), restricted to $I(\chi)^B$ ($B$-invariants), then these
form a basis for the dual of $I(\chi)^B$. This was proven by
Casselman in \cite{C}.

Exploiting this fact involves yet another complication: $I(\chi)$
has two natural bases: One is $\{\phi_{BwB}\}_w$, which has
already been introduced. This basis is suitable for computations
using integral expressions like (\ref{integral}). The second one
is the basis which is dual to the functionals coming from
intertwining operators that were mentioned above. This is the
``Casselman basis'' $\{f_w\}_w$. This basis is useful if we have
already computed the effect of the Shalika functional on an
element of this basis and wish to extend the computation to all
elements.

It is essential to establish a connection between the two bases.
The only immediate relation is that $\phi_{Bw_l B}=f_{w_l}$.
Therefore, a good starting point would be to compute the effect of
the Shalika functional on $g_\lambda$-translates of $\phi_{Bw_l
B}$. This is the approach originally followed by Casselman and
Shalika when computing Whittaker vectors.

However, for most of the subgroups $H$ that we are interested in, we cannot
expect $B w_l B$ - translates to belong to a single double $P
\backslash G / H$ coset, either. It will, on the contrary, be
usually the case (and the original work of Casselman and Shalika,
as well as much of similar subsequent work, can be reformulated in
these terms) that suitable translates of $B$ \emph{will} belong to
a single double coset. This makes it possible to compute the
effect of the Shalika functional on translates of $\phi_B$. But
now we have the problem of connecting to the Casselman basis. The
work of Hironaka \cite{H} shows how to do that.

The basic new idea is that, instead of expressing $\phi_K$ in the
Casselman basis, one expresses \emph{the projection of the
distribution $\Delta_H$ to $B$-invariants} in that basis. We explain it below:

There is a natural projection from the space of distributions $D$
satisfying (\ref{distr}) to $I(\chi^{-1})^B$ given by:
\begin{equation} \label{compute}
 R_B D (g)= \int_B D(gb) db = \int_G
 D(b) 1_B (g^{-1}b) db.
 \end{equation}

It will not impede our computation to apply this projection, since
\begin{equation}\Omega_H(g)= \Delta_H ( R_g 1_K ) =
R_{g^{-1}}\Delta_H(1_K) =
R_B R_{g^{-1}}\Delta_H(1_K)
\end{equation}
(here we have used the fact that $R_{g^{-1}}=R_g^*$, the adjoint
of $R_g$, that $R_B^*=R_B$ and that $1_K$ is $B$-invariant)

Our goal will be to express $R_B R_{g^{-1}}\Delta_H$ in the
Casselman basis:
\begin{equation}\label{expression}
   R_B R_{g^{-1}} \Delta_H = \sum_w a_w ( g ) f_w.
\end{equation}

Based on the work of Casselman and Hironaka (we refer the reader
to \cite{H} or \cite{O} for details), the coefficients $a_w$ will
be given by:
\begin{equation} \label{coeffs}
 a_w ( g ) = \frac{c_w ( \chi^{- 1} )}{c_{w^{- 1}} (^w
\chi )}
     T_{w^{- 1}}^{\ast} \Delta_H ( R_{g} 1_{B, ^w\chi} )
\end{equation}
where
\begin{equation}
c_w(\chi)=\prod_{\alpha>0,w\alpha<0}c_\alpha(\chi)
\end{equation}
and
\begin{equation} \label{cadefinition}
c_\alpha(\chi)=\frac{1-q^{-1}\chi(a_\alpha)}{1-\chi(a_\alpha)}
=\frac{1-q^{-1}e^\alpha}{1-e^\alpha}(g_\chi)
\end{equation}
so we need to compute $T_{w^{- 1}}^{\ast} \Delta_H ( R_{g} 1_{B,
^w\chi} )$. Notice that by the correspondence between functionals on $I(\chi)$ and certain distributions on $G$, we freely apply $T_w^*$, the adjoint of $T_w$, to the distribution $\Delta_H$.

Finally, Casselman computed in \cite{C},\textsection 4 the effect
of $f_{w,\chi^{-1}}$ (now thought of as an element of $I(\chi)^*$)
on $\phi_K$:
\begin{equation}\label{fwK}
f_{w,\chi^{-1}}(\phi_{K,\chi})=Q^{-1}\frac{c_{w_l}(^w\chi)}{c_w(\chi^{-1})}
\end{equation}
where $Q$ is some constant independent of $\chi$. More precisely, $Q$ is the harmonic mean of the numbers $(BwB:B)\, , w\in W$, i.e.:
\[Q^{-1}=\sum_w (BwB:B)^{-1} = \frac{\tmop{meas}(Bw_l B)}{\tmop{meas(K)}}.\]

Therefore, knowledge of the coefficients $a_w(g)$ allows us to
compute $\Omega(g)$. If we combine the equations above, the $c$-factors from (\ref{goal}) and (\ref{fwK}) will simplify to
give:
\[\frac{c_{w_l}(^w \chi)}{c_{w^{-1}}(^w \chi)}\]
which equals
\[\prod_{\alpha\in \Phi_{GL}^+, w\alpha>0} c_\alpha(\chi)\]
so finally we get:

\begin{equation} \label{goal}
\Omega_H(g)=Q^{-1}\sum_w\prod_{\alpha\in \Phi_{GL}^+, w\alpha>0} c_\alpha(\chi)
T_{w^{-1}}^*\Delta_H(R_g 1_{B,{^w\chi}}).
\end{equation}


\section{Vanishing results and the integral on the open orbit}
\label{secvanishing}

Remember that
\[ \Omega(g_\lambda)= \Omega_H \left( \left(
\begin{array} {cc} \varpi^{-\lambda} & \\ & I \end{array} \right )
\right)\] and $\lambda_1\ge\lambda_2\ge\dots\ge\lambda_n$.

The role of the character $\Psi$ is to make $\Omega(g_\lambda)$
vanish if not all $\lambda_i\ge 0$. It will appear soon why this
is crucial for our method. Let $X\in \tmop{Mat}_n(\mathfrak o)$.
We have:

\[\Omega(g_\lambda)=\Omega \left(g_\lambda \cdot \left(
\begin{array} {cc}  I & X\\ & I \end{array} \right )\right)
 = \Omega \left( \left(
\begin{array} {cc} I & \varpi^\lambda X\\ & I \end{array} \right )
\cdot g_\lambda \right) =\psi(\tmop{tr}(\varpi^\lambda X))\Omega
(g_\lambda)\] so if $\varpi^\lambda \notin \tmop{Mat}_n (\mathfrak
o)$ we can find $X$ such that $\psi(\tmop{tr}(\varpi^\lambda X))
\ne 0$. From this it follows that $\Omega(g_\lambda)=0$.

We are left with computing $T_{w^{- 1}}^{\ast} \Delta_H ( R_{g}
1_{B, ^w\chi} )$ for $g=  \left(
\begin{array} {cc} \varpi^{-\lambda} & \\ & I \end{array}
\right)$ where $\lambda_1\ge \dots \ge \lambda_n \ge 0$. The
function $R_g 1_B$ is supported on the set $Bg^{-1}$.

\begin{lemma} \label{basiclemma}
For those $g$, $Bg^{-1}\subset P \cdot H$.
\end{lemma}

\begin{proof}
By the Iwahori factorization for $B$, $B=P_0 N_1^-$, it suffices
to show that $N_1^- g^{-1} \subset P \cdot H$. But for $g$ as
above, $gN_1^-g^{-1}\subset N_1^{-1}$ and $g \in P$, therefore
there remains to show that
\begin{equation} \label{N1}
N_1^-\subset P \cdot H.
\end{equation}

We show something stronger, because it will be needed later.
Namely, we prove that
\begin{equation}\label{n1}
N_1^{-\alpha}\subset P_0^{\alpha,\tilde\alpha} H_0
\end{equation}
where $P_0^{\alpha,\tilde\alpha}=P_0^{\alpha}\cdot P_0^{\tilde\alpha}$.

This is essentially a simple approximation argument on the Lie
algebra: Write a given element of $N_1^-$ as $I+n_1$, then the
entries of $n_1$ will be in $\mathfrak{p}$. We argued above that
the Lie algebras of $H$ and $P$ are complementary, therefore we
can find an integral matrix $h_1$ of the form (\ref{liealgebra}),
with coefficients in $\mathfrak{p}$ such that $p_1=n_1+h_1$ is
upper triangular. In fact, we can arrange so the entries of $p_1$
above the diagonal will only be non-zero in the $\alpha$ and
$\tilde\alpha$ positions. We will then have
$(I+n_1)(I+h_1)=I+n_1+h_1+n_1 h_1=$upper triangular+$n_2$, where
$n_2$ has coefficients in $\mathfrak{p}^2$. Then similarly we will
find $h_2$ with coefficients in $\mathfrak{p}^2$ of the form
(\ref{liealgebra}) such that $n_2+h_2$ is upper triangular, etc,
and then the converging sum $h=I+h_1+h_2+\dots$ will satisfy:
\[(I+n_1)h\in P_0^{\alpha,\tilde\alpha} \]
\end{proof}

For elements of $I(\chi)$ with support in $P\cdot H$, $\Lambda_H$
will have a very simple form. We first prove a vanishing result
as in \cite{MR} and \cite{O}:

\begin{proposition} \label{propvanishing}
For $w \notin \Gamma$, the distribution $T_{w^{-1}}^* \Delta_H$ is
supported away from $P\cdot H$.
\end{proposition}

Remember that $\Gamma$ denotes the Weyl group of $\tmop{Sp}_{2n}$,
considered as a subgroup of $W$.

\begin{proof}
The argument is exactly that of the aforementioned papers:
Computing formally at first (treating the distribution as a
function), we know that $T_{w^{-1}}^* \Delta_H$ satisfies:
\[ T_{w^{-1}}^* \Delta_H ( p g h ) = {^w\chi}^{- 1}
\delta^{\frac{1}{2}} ( p ) T_{w}^* \Delta_H ( g ) \Psi_H^{- 1} ( h
) \] for $p \in P, h \in H$. Hence, if $x\in P \cap H$ we get:
\[ \Psi_H^{-1}(x) T_{w^{-1}}^* \Delta_H (1) = T_{w^{-1}}^* \Delta_H (x)
= {^w\chi}^{-1} \delta ^{ \frac{1}{2}} (x) T_{w^{-1}}^* \Delta_H
(1).\]

The group $P \cap H$ consists of the matrices of the form:
\begin{equation} \label{PcapH}
 x=\left(\begin{array}{cccccc}
\alpha_1\\
&\ddots\\
&&\alpha_n \\
&&& \alpha_n\\
&&& &\ddots&\\
&&& &&\alpha_1 \end{array}\right).
\end{equation}

Therefore $\Psi_H(x)=1$ for all $x$, while the character $^w
\chi^{- 1} \delta^{\frac{1}{2}}$ will be trivial on all such $x$
if and only if $w \in \Gamma$.

To make this rigorous, consider the space of all $f\in I(^w\chi)$
that are supported in $P \cdot H $. Restriction to $H$ provides an
isomorphism between the space of such $f$ and $\mathcal{S}(P \cap
H\backslash H,r)$, the Schwartz-Bruhat space of smooth functions on $H$
which are compactly supported modulo $P \cap H$ and vary on left
multiplication by $P \cap H$ via the twisting character
$r={^w\chi}\delta^{\frac{1}{2}}$. This character is trivial if and
only if $w\in\Gamma$.

The functional $T_{w^{-1}}^*\Lambda_H$, restricted to this space,
can be lifted to a distribution $D$ on $H$ via the H-equivariant
projection: $\mathcal P_r: C_c^\infty(H) \to \mathcal{S}(P \cap H
\backslash H,r)$ where \[\mathcal P_r (f)(h)=\int_{P \cap H}
r^{-1}(h_0)f(h_0h)dh_0.\] Notice that $P\cap H$ is unimodular.

$D$ has the property: $D (R_h f)=\Psi_H(h) D(f)$. It follows that
$\Psi_H \cdot D$ is a right-invariant distribution on $H$, and
since $H$ is unimodular it is both left and right Haar measure. In
other words:

\begin{equation} \label{D} D = \Psi_H^{-1} dh. \end{equation}

Using the above relation, for $x\in P \cap H$, $f\in
C_c^\infty(H)$
\begin{eqnarray}
 \int_H f(xh) D(h) dh = \int_H f(xh) \Psi_H^{-1} (h)
dh = \nonumber\\ = \Psi_H(x) \int_H f(h) \Psi_H^{-1}(h)dh = \int_H
f(h) D(h) dh.
\end{eqnarray}
On the other hand, $D$ is supposed to factor through
$\mathcal{S}(P \cap H\backslash H,r)$, and we have $\mathcal P_r
(f(x\cdot))= r(x) \mathcal P_r (f(\cdot))$, therefore $D$ has to
be zero unless $r=1$.

\end{proof}

\begin{corollary}
For $w\notin \Gamma$, $a_w (g_\lambda)=0$ for every $\lambda$.
\end{corollary}
\begin{proof}
It follows immediately from the above lemma and (\ref{goal}).
\end{proof}

\begin{corollary}
If $\chi$ is not of the form (\ref{chi}) (or a $W$-translate of
this) then there exists no Shalika model for $I(\chi)$.
\end{corollary}
\begin{proof}
Indeed, following the proof of the proposition, all Casselman
coefficients in that case would vanish.
\end{proof}

\begin{corollary}
For $w\in \Gamma$, $\phi \in I({^w\chi})$ with $\tmop{supp}\phi \subset
P\cdot H$,
\begin{equation} \label{Lambda}
\Lambda_H(\phi)= \int_{P \cap H \backslash H} \phi(h)
\Psi_H^{-1}(h) dh
\end{equation}
\end{corollary}
\begin{proof}
This follows from (\ref{D}).
\end{proof}
Notice that the functional defined by this integral is clearly non-zero, as it lifts to the non-zero distribution $\Psi_H^{-1}dh$ on $C_c^\infty(H)$.


\section{Dependence on $\lambda$} \label{secdependence}

There is an alternative expression to (\ref{Lambda}) which is
going to be useful later: If $\phi=\mathcal P_\chi (f)$ with $f
\in C_c^\infty(G)$ then combining the integral expression
(\ref{Pchi}) for $\mathcal P_\chi$ with that of (\ref{Lambda}) we
get:
\begin{equation} \label{Delta}
\Delta_H(f)= \int_{PH} \chi^{-1}\delta^{\frac{1}{2}} \left ( p(x)
\right) f(x) \Psi_H^{-1} \left( h(x) \right) dx.
\end{equation}
Here the measure is Haar measure on $G$; remember that $P\cdot H$
is open and dense in $G$. The symbols $p(x)\in P$ and $h(x)\in H$ correspond
to a factorization of $x$: $x=p(x)h(x)$. They are only well
defined modulo $P \cap H$ but that doesn't matter since the
characters are trivial there.

Recall that up to this point we have not specified a
normalization for the Shalika functional, since the expressions we
have considered are only determined up to a constant. \emph{We now
fix a Haar measure on $G$ such that the measure of the Iwahori
subgroup $B$ is 1. The normalization for the Shalika functional
will then be that corresponding to (\ref{Delta}).} As we shall see
immediately, this normalization corresponds to
\begin{equation}\label{normalization}
\Lambda(\phi_B)=1.
\end{equation}
Note that this will be our working convention, but for the sake of
a simpler expression \emph{the normalization changes when we state
our main theorem}.

\begin{proposition} \label{depproposition}
For $g=g_{-\lambda}=\left(\begin{array} {cc} \varpi^{-\lambda} \\
&I
\end{array} \right)$, $\lambda_1\ge \dots \ge \lambda_n\ge 0$, we have
$\Lambda_H(R_g \phi_B)= \chi^{-1}\delta^{\frac{1}{2}}(g_\lambda)$.
\end{proposition}

\begin{proof}
Use the Iwahori factorization to write an arbitrary element $x\in
Bg^{-1}$ as $x=p_0 n_1^- g^{-1}= \left( p_0 g^{-1} \right) \left(
g n_1^- g^{-1} \right)$. Since $g n_1^- g^{-1} \in N_1^-\subset
P_0 H_0$ we get that $\Psi_H\left(h(x)\right)=1$ and
$\chi^{-1}\delta^{\frac{1}{2}} \left ( p(x) \right) =
\chi^{-1}\delta^{\frac{1}{2}} (g^{-1}) =
\chi^{-1}\delta^{\frac{1}{2}} (g_\lambda)$. Therefore, using
(\ref{Delta}) we have
\[\Delta_H(R_g 1_B) = \chi^{-1}\delta^{\frac{1}{2}} (g_\lambda)
\cdot \tmop{Vol}(Bg^{-1}) = \chi^{-1}\delta^{\frac{1}{2}}
(g_\lambda).\]
\end{proof}

This gives the coefficients $a_1(g_{-\lambda})$ by (\ref{coeffs}). It
actually gives more: Since $T_{w^{-1}}$ is $G$-equivariant and the
Shalika functional is unique up to scaling, $T_{w^{-1}}^*
(\Lambda)$ will be a multiple of $\Lambda_{^w\chi}$ (normalized as
in (\ref{normalization})) therefore by means of (\ref{coeffs}) the
above considerations prove:

\begin{corollary} For every $w \in \Gamma$, $\lambda$ as above,
\begin{equation} \label{dependence}
\frac{a_w(g_{-\lambda})}{a_w(g_0)} =
{^w}\chi^{-1}\delta^{\frac{1}{2}} (g_\lambda).
\end{equation}
\end{corollary}

Applying this to (\ref{goal}) we get:
\begin{equation}\label{general}
\Omega_H(g_{-\lambda})=Q^{-1}\sum_w\prod_{\alpha\in \Phi_{GL}^+, w\alpha>0} c_\alpha(\chi) {^w\chi}^{-1}\delta^{\frac{1}{2}}(g_\lambda)
T_{w^{-1}}^*\Delta_H(1_{B,{^w\chi}}).
\end{equation}


\section{Analytic results} \label{secanalytic}

In this section we establish two important analytic results:
First, the convergence of the period
integral (\ref{Lambda}) for \emph{all} $\phi\in I(\chi)$, in the
case that the induction data lie in a certain open region. Second, the rationality of the Shalika function with respect to
the Satake parameters $(q^{\pm z_1},\dots,q^{\pm z_n})$,
which will allow us to restrict our attention to the region of convergence.

\begin{proposition} \label{period}
When $\tmop{Re}z_1>\tmop{Re}z_2>\dots>\tmop{Re}z_n$, the period
integral (\ref{Lambda}) converges absolutely for every $\phi \in
I(\chi)$, and therefore represents a (non-zero) Shalika functional.
\end{proposition}

\begin{proof}
This is the only point where it will be more convenient to refer
to the Shalika subgroup itself, rather than $H$.

Using the correspondence of (\ref{translate}), the equivalent to
(\ref{Lambda}) integral for $\Lambda$ is:
\begin{eqnarray*}& \Lambda(\phi)= \int_{_P S \backslash S} \phi(\xi s)
\Psi^{-1}(s)
ds = \nonumber \\ & \int_{T \backslash \tmop{GL}_n(F)}
\int_{\tmop{Mat}_n(F)} \phi \left(\left( \begin{array}{cc} &I \\
w_0\end{array}\right) \left( \begin{array}{cc} I& X\\ &
I\end{array}\right) \left(
\begin{array}{cc} g \\ &g\end{array}\right)\right)\psi^{-1}(\tmop{tr}X)
dX \,\,
dg.
\end{eqnarray*}
Here $_P S= \xi^{-1} P \xi \cap S$, and $T$ is the maximal torus
of diagonal matrices in $\tmop{GL}_n(F)$.

Since every element of $I(\chi)$ is a locally constant function,
which is determined by its restriction to $K$ and therefore
dominated by a suitable multiple of $\phi_K$, it suffices to prove
the proposition for $\phi=\phi_K$. Using an Iwasawa decomposition
for $\tmop{GL_n}(F)$, we can write $T \backslash \tmop{GL}_n(F) =
U \cdot K_n$ (measure-preserving) where $U$ is the subgroup of
upper triangular $n\times n$ matrices with 1's on the diagonal and
$K_n=\tmop{GL}_n(\mathfrak o)$. Given the $K$-invariance of
$\phi_K$, the above integral reduces to:

\begin{multline}   \Lambda(\phi_K)=  \nonumber \\ \int_U
\int_{\tmop{Mat}_n(F)} \phi_K \left(\left( \begin{array}{cc} &I\\
w_0\end{array}\right) \left(
\begin{array}{cc} I& X\\ & I\end{array}\right) \left(
\begin{array}{cc} n \\ &n\end{array}\right)\right) \psi^{-1}(\tmop{tr}X)
dX \,\,
dn \nonumber \\
= \int_U\int_{\tmop{Mat}_n(F)}  \phi_K \left(\left( \begin{array}{cc} I&\\
&n\end{array}\right)  \left(\begin{array}{cc} &I\\w_0
\end{array}\right) \left(
\begin{array}{cc} n \\ &I\end{array}\right)\left(
\begin{array}{cc} I& n^{-1}Xn\\ & I\end{array}\right)\right)\nonumber \\
\psi^{-1}(\tmop{tr}X) dX \,\, dn.
\end{multline}

The factor $\left( \begin{array}{cc} I&\\
&n\end{array}\right)$ on the left disappears because $\phi_K \in
I(\chi)$, and $n^{-1}Xn$ can be replaced by $X$ since conjugation
by $n$ is a measure preserving automorphism of $\tmop{Mat}_n(F)$.
Therefore, the integral above is dominated absolutely by the
integral which represents the intertwining operator for the Weyl
group element $w=\left( \begin{array}{cc} &w_0\\
I&\end{array}\right)$:
\begin{equation} \label{intw}
 T_w(\phi) = \int_{U'} \phi(w^{-1}u) du
\end{equation}
where $U'$ is the group of upper triangular unipotent matrices
with the identity element in the lower $n\times n$ block.

It is known that the integral (\ref{intw}) converges absolutely
for $\tmop{Re}z_1>\tmop{Re}z_2>\dots>\tmop{Re}z_n$, which
establishes the claim.
\end{proof}

Now, let $D$ denote the algebraic variety of diagonal elements
in $\tmop{GL}_n(\mathbb C)$, identified as above with the set of
unramified characters $\chi$ of the form (\ref{chi}). $\mathbb C[D]$
will denote the algebra of regular functions on $D$ and $L$ its quotient
field. Let $X$ be the space of all locally
constant functions on $K$, which are left invariant under $P_0$.
For every $\chi$, $X$ can be identified with $I(\chi)$ via
restriction of functions in $I(\chi)$ to $K$. It then makes sense to talk
about a rational family $\{\phi_\chi\}_{\chi\in D}$ with $\phi_\chi \in
I(\chi)$, in the sense that $\phi \in
L\otimes_\mathbb C X$.

\begin{proposition} \label{holomorphicity}
There exists a non-zero Shalika functional for almost all $\chi\in D$. Moreover, if $\{\phi_\chi\}_{\chi\in D}, \,\phi_\chi\in I(\chi)$ is a rational
family then $\Lambda_{H,\chi}(\phi_\chi)$ is
a rational function of $\chi\in D$.
\end{proposition}

This will be a direct application of a theorem of Bernstein. We
simply state Bernstein's theorem, and refer the reader to
\cite{Be} or \cite{GPS}, p.127, for explanations and the proof:

\begin{theorem*}[Bernstein]
Let $X$ be a vector space of countable dimension over $\mathbb C$,
$X^*$ its linear dual, $D$ an irreducible variety over $\mathbb
C$, $\mathbb{C}[D]$ the algebra of regular functions on $D$, $L$
its quotient field, $X_L=L\otimes_{\mathbb C}X$ and
$X_L^*=\tmop{Hom}_L(X_L,L)$. For every $d\in D$, consider a system of
linear equations on $X^*$:
\begin{equation}\label{system}
\{\left<x_{rd},\Lambda\right>=l_{rd}\}_{r\in R}\,\, ,\,\, x_{rd}\in X,
l_{rd}\in \mathbb C
\end{equation}
($R$ a fixed indexing set for the equations of the system), and assume that $x_r:=\{x_{rd}\}_d$ and $l_r:=\{l_{rd}\}_d$ vary rationally with $d$, i.e.
$x_r\in X_L$ and $l_r \in L$. Assume that for $d$ in
some $\Omega\subset D$, open in the usual topology, the system
(\ref{system}) has a unique solution $\Lambda_d\in X^*$. Then the
system (\ref{system}), considered as a system of linear equations on $X_L^*$,  has a unique solution $\Lambda\in X_L^*$.
\end{theorem*}

To prove the proposition, consider the system of equations which consists
of the requirements:

-$\Lambda$ is a Shalika functional, i.e. equation (\ref{Shalika})
for all $f\in X$

-normalization condition: equation (\ref{normalization}).

It is easy to see that these equations are rational in $\chi$. The
solution is
then a (normalized) Shalika functional. Proposition \ref{period} shows that a Shalika function exists in the case that the
inducing parameters lie in the aforementioned region, and we also know that this functional is unique. By Bernstein's theorem, $\Lambda$ extends to an element of $X_L^*$. Therefore, when
applied to $\phi_\chi$, the result will be a rational function in $\chi$,
in other words a rational function in $q^{z_1},\dots,q^{z_n}$.

The phrase ``almost all $\chi$'' refers, of course, to the possible singular hypersurfaces of $\Lambda$. At the end of our proof we will become more precise about where these might lie.


\section{The functional equations}\label{secfe}

We have already established the dependence of the Casselman
coefficients on $\lambda$. There remains to determine the
dependence on $w$. By the uniqueness of Shalika models, we know
that $T_{w^{-1}}^*\Delta_{H,\chi}$ has to be a constant multiple
of $\Delta_{H,^w\chi}$. We will use the integral expression
provided by Proposition \ref{period} to compute this constant
explicitly for $w$ in a set of generators in $\Gamma$. That set of
generators consists of the transposition $(n, n+1)$ (which is a
simple reflection corresponding to a long root of
$\tmop{Sp}_{2n}$) and of the elements $(i, i+1)(2n-i, 2n+1-i)$
with $1\le i <n$ (which are simple reflections in the Weyl group
of $\tmop{Sp}_{2n}$ corresponding to short roots). By expressing
an arbitrary element of $\Gamma$ as a product of simple
reflections, and writing the intertwining operators as a
composition of intertwining operators correspondingly, the result
will follow for all $w\in \Gamma$. To be rigorous, we cannot
iterate the explicit computation, since after applying the first
intertwining operator the inducing parameters will no longer
belong to the region of convergence for the period integral.
However, rationality will allow us to extend the results to the
region of non-convergence.

(Notice that no functional equations for $\Omega(g)$ will appear
explicitly. The title of this section serves as a connection to
the Casselman-Shalika method as used in the literature.)

To compute $T_{w^{-1}}^*\Lambda_{H,\chi}$ as a multiple of
$\Lambda_{H,^w\chi}$ it suffices to compute their quotient when
applied to a single element, for instance $\phi_B$. Here we use
Theorem 3.4 of \cite{C}, to write for any simple reflection in $W$
(corresponding to the root $\alpha$):
\begin{equation}\label{simple}
T_{w_{\alpha}}(\phi_B)=(c_\alpha (^{w_\alpha}\chi)-1)\phi_B +
q^{-1}\phi_{Bw_a B}
\end{equation}
and for a simple $\tmop{Sp}_{2n}$-reflection of a short root
$w=w_\alpha w_{\tilde{\alpha}}$:
\begin{eqnarray}\label{double}
& T_w(\phi_B)=
(c_\alpha(^{w}\chi)-1)(c_{\tilde{\alpha}}(^{w}\chi)-1)\phi_B +
q^{-1}(c_\alpha(^{w}\chi)-1)\phi_{Bw_{\tilde{\alpha}}B} +
\nonumber
\\& +q^{-1}(c_{\tilde{\alpha}}(^{w}\chi)-1)\phi_{Bw_{\alpha}B} +
q^{-2} \phi_{BwB}.
\end{eqnarray}

Now we need to apply $\Lambda_H$ to these expressions to get the
functional equations.

\begin{proposition}
Let $w=w_\alpha$ be the simple reflection $(n, n+1)$ in $\Gamma$.
Then
\begin{equation}\Lambda_{H,\chi}(T_{w^{-1}}\phi_{B,^w\chi})=(-1)\chi(a_\alpha)c_\alpha(\chi).
\end{equation}
\end{proposition}

\begin{proof}
By (\ref{simple}), we need to compute $\Lambda_H(\phi_{Bw_\alpha
B})$. Assume that $\tmop{Re}z_1>\tmop{Re}z_2>\dots>\tmop{Re}z_n$,
so that $\Lambda_H$ is given by (\ref{Lambda}). In order to apply
(\ref{Lambda}) or (\ref{Delta}) we need to express a generic
element of $B w_{\alpha} B$ in the form $P\cdot H$. (Remember that
the latter is open and dense, and since $B w_{\alpha} B$ is open,
almost every element will be expressible in this form.) We use the
factorization
\[N=\prod N^\beta\]
with $\beta$ running over all positive roots (in any order), and
the similar factorizations for $N_0$, $N_1^-$, in order to write
the double coset $Bw_\alpha B$ in a measure-preserving way (except
for a constant factor $[Bw_\alpha B : B]= q$) as:
\begin{equation}
P_0 w_\alpha N_0^{\alpha} N_1^{-,\widehat\alpha}
\end{equation}
where all compact groups that appear are assumed to have measure 1
and $N_1^{-,\widehat\alpha}$ denotes that the factor corresponding
to $\alpha$ is missing.

Remember that by (\ref{n1}),
\begin{equation}
N_1^{-,\widehat\alpha}\subset P_0^{\widehat\alpha}H_0
\end{equation}
in other words the $P_0$ factor contains no $N^\alpha$ factor.
Therefore it can be pulled to the other side of $w$ to produce a
$P_0$ factor on the left. We therefore have:
\begin{equation}
\Delta_H(1_{B w_\alpha B})= q \int_{w_\alpha N_0^\alpha}
\chi^{-1}\delta^{\frac{1}{2}} \left ( p(x) \right) \Psi_H^{-1}
\left( h(x) \right) dx
\end{equation}
(where of course the integrand is only defined on a dense open
subset of $w_\alpha N_0^\alpha$).

For an element $x\in w_\alpha N_0^\alpha$ we compute the
factorization:
\begin{eqnarray}\label{Nalpha} \scriptstyle
& x = \left(\begin{array}{cccccc} 1\\ &\ddots \\ &&0 & 1 \\  &&1&y \\
&&&& \ddots \\ &&&&& 1 \end{array} \right) = \nonumber \\ &
\left(\begin{array}{cccccc} 1\\ &\ddots \\ &&-y^{-1} & 1 \\  &&&y \\
&&&& \ddots \\ &&&&& 1 \end{array} \right) \cdot
\left(\begin{array}{cccccc} 1\\ &\ddots \\ &&1 &  \\  &&y^{-1}&1 \\
&&&& \ddots \\ &&&&& 1 \end{array} \right) \in P\cdot H \nonumber
\\
\end{eqnarray}
if $y\ne 0$.

{}From this we see that $\chi^{-1}\delta^{\frac{1}{2}} (p(x))=
|y|^{2z_n-1} =\left(q\chi (a_\alpha)\right) ^{\tmop{val}(y)}$ and
\\ $\Psi_H^{-1} (h(x)) = \psi^{-1}(y^{-1})$, therefore
\begin{eqnarray}
\Delta_H(1_{B w_\alpha B})= q \int_{\mathfrak o} |y|^{2z_n-1}
\psi^{-1}(y^{-1}) dy = \nonumber \\
q \sum_{i=0}^\infty \left(q\chi(a_\alpha)\right)^i
\int_{\mathfrak{p}^i - \mathfrak{p}^{i+1}} \psi^{-1}(y^{-1}) dy.
\end{eqnarray}

Using the fact that the conductor of $\psi^{-1}$ is $\mathfrak o$,
we have:
\begin{equation}
\int_{\mathfrak{p}^i} \psi^{-1} (y) dy = \left\{\begin{array}{ll}
\tmop{Vol}(\mathfrak{p}^i)&\textrm{, if  } i\ge 0 \\
0 & \textrm{, otherwise} \end{array}\right.
\end{equation}

Thus, by making a change of variables in the above expression, all
integrals vanish except for the integral on
$\mathfrak{p}^0-\mathfrak{p}^1=\mathfrak{o}^\times$, which is
equal to $(1-q^{-1})$, and the integral on
$\mathfrak{p}^1-\mathfrak{p}^2$, which is equal to $q^{-2}\cdot
(-1)$. Hence finally:
\begin{equation}
\Delta_H(1_{Bw_\alpha B}) = q \cdot
(1-q^{-1}-q^{-1}\chi(a_\alpha)).
\end{equation}

Therefore, applying $\Lambda_H$ to equation (\ref{simple}) and
after some simple algebraic manipulation we get:
\begin{equation}
T_{w^{-1}}^*\Lambda_{H,\chi}(\phi_{B,^w \chi}) = (-1)
\chi(a_\alpha) c_\alpha(\chi).
\end{equation}

This establishes the result for $\chi$ in a certain region. As
shown in \cite{Be} (essentially in the same way that we proved
proposition \ref{holomorphicity}), the intertwining operators are
rational in $\chi$, hence $T_{w^{-1}}\phi_{B,^w\chi}$ is a
rational family (in the sense of proposition
\ref{holomorphicity}) so by proposition \ref{holomorphicity} the
result follows for all $\chi$.

\end{proof}

Now for the simple reflections corresponding to short roots:

\begin{proposition}
Let $w=w_\alpha w_{\tilde{\alpha}}$, where $w_\alpha=(i, i+1)$ ,
$1\le i <n$. Then
\begin{equation}\Lambda_{H,\chi}(T_{w^{-1}}\phi_{B,^w
\chi})=(-1)\chi(a_\alpha)
\frac{1-q^{-1}\chi(a_{-\alpha})}{1-q^{-1}\chi(a_\alpha)}
c_\alpha(\chi)c_{\tilde{\alpha}}(\chi).\end{equation}
\end{proposition}

\begin{proof}
Assume again that $\tmop{Re}z_1>\tmop{Re}z_2>\dots>\tmop{Re}z_n$.
We are going to apply $\Lambda_H$ to (\ref{double}), but first we
use a trick to reduce the number of computations needed:

Recall that $T_{w_\alpha^{-1}}^* \Delta_H$ is supported away from
$P\cdot H$, while the support of $\phi_{B,^w\chi}$ is contained in
$P\cdot H$. Therefore $T_{w_\alpha^{-1}}^* \Lambda_H
(\phi_{B,^w\chi})=0$. Substituting from (\ref{simple}) yields:
\begin{equation}
\Lambda_{H,\chi}(\phi_{Bw_\alpha B, \chi})=q(1-c_\alpha(^{w_\alpha}
\chi))
\end{equation}
and similarly for $\phi_{Bw_{\tilde{\alpha}}B,\chi}$ (just replace
$\alpha$ by $\tilde{\alpha}$).

Substituting in (\ref{double}) and simplifying, we get:
\begin{equation} \label{doubledouble}
T_{w^{-1}}^*\Lambda_{H,\chi}(\phi_{B,^w \chi}) =
-(1-c_\alpha(^{w_\alpha}\chi))^2 + q^{-2} \Lambda_{H,\chi}
(\phi_{BwB,\chi}).
\end{equation}

Therefore we only need to compute $\Lambda_{H,\chi}
(\phi_{BwB,\chi})$.

As in the previous proposition, we write in a measure preserving
way (except for a constant factor $[BwB:B]=q^2$):
\begin{equation}
BwB=P_0 w N_0^\alpha N_0^{\tilde{\alpha}}
N_1^{-,\widehat\alpha,\widehat{\tilde{\alpha}}}.
\end{equation}

The $N_1^{-,\widehat\alpha,\widehat{\tilde{\alpha}}}$ factor
belongs to $P_0^{\widehat\alpha,\widehat{\tilde{\alpha}}}\cdot H$,
and the $P_0^{\widehat{\alpha},\widehat{\tilde{\alpha}}}$ factor
can be pulled to the left without changing the measure on
$N_0^\alpha N_0^{\tilde{\alpha}}$. Finally, using the fact that
the permutation matrix representing $w$ belongs itself in $H_0$,
we have:
\begin{equation}
P_0 w N_0^\alpha N_0^{\tilde{\alpha}}
N_1^{-,\widehat\alpha,\widehat{\tilde{\alpha}}} \subset P_0 w
N_0^\alpha N_0^{\tilde{\alpha}} w^{-1} H_0.
\end{equation}

A generic element in $w N_0^\alpha N_0^{\tilde{\alpha}} w^{-1}$
can be factored in the form $P\cdot H$ as follows:
\begin{eqnarray} &
x = \left(\begin{array}{ccccccc} \ddots \\ &1&\\ &y_1&1 \\
&&&\ddots\\  &&&&1 \\ &&&& y_2 &1 \\ &&&&&& \ddots \end{array}
\right) = \nonumber
\end{eqnarray}
\begin{eqnarray} & = \left(\begin{array}{ccccccc}
\scriptstyle \ddots \\
& \scriptstyle (1-y_1 y_2)^{-1}
& \scriptstyle -y_2(1-y_1 y_2)^{-1}\\ && \scriptstyle 1 \\
&&& \scriptstyle \ddots\\ &&&&
\scriptstyle (1-y_1 y_2)^{-1} & \scriptstyle -y_1(1-y_1 y_2)^{-1}\\
&&&&& \scriptstyle 1 \\
&&&&&& \scriptstyle \ddots
\end{array} \right)\nonumber \\ &\cdot
\left(\begin{array}{ccccccc} \ddots \\ &1&y_2\\ &y_1&1 \\
&&&\ddots\\ &&&&1 & y_1 \\ &&&& y_2 &1 \\ &&&&&& \ddots
\end{array} \right) \in P\cdot H
\end{eqnarray}
(if $y_1 y_2 \ne 1$).

{}From this we see that $\chi^{-1}\delta^{\frac{1}{2}} (p(x))=
|1-y_1 y_2|^{z_i-z_{i+1}-1}$ and $\Psi_H^{-1} (h(x)) = 1$,
therefore
\begin{eqnarray}&
\Delta_H(1_{B w_\alpha B})= q^2 \int_{\mathfrak{o\times o}} |1-y_1
y_2|^{z_i-z_{i+1}-1} dy_1 dy_2 = \nonumber \\& =q^2 \left[
\int_{y_1 \in \mathfrak p} \int_{y_2 \in \mathfrak o} 1 dy_1 dy_2
+ \int_{y_1 \in \mathfrak o^\times} \int_{y_2 \in \mathfrak p} 1
dy_1 dy_2 + \right.\nonumber \\ &\left.+\int_{y_1, y_2 \in
\mathfrak o^\times} |1-y_1y_2|^{z_i-z_{i+1}-1} dy_1dy_2 \right]=
\nonumber
\\& = q^2 \left[ q^{-1}+q^{-1}(1-q^{-1}) +
(1-q^{-1})\int_{\mathfrak o^\times} |1-u|^{z_i-z_{i+1}-1} du
\right].
\end{eqnarray}

We split the last integral into $1-u \notin \mathfrak p$ and $1-u
\in \mathfrak p$. The former contributes $\frac{q-2}{q}$. For the
latter, we substitute $t=1-u$ and integrate over $\mathfrak p$ to
get:
\begin{equation}
\sum_{j=1}^\infty q^{-j(z_i-z_{i+1})}(1-q^{-1}) = (1-q^{-1})
\frac{q^{-(z_i-z_{i+1})}}{1-q^{-(z_i-z_{i+1})}}.
\end{equation}

Putting together all the above we eventually get what the
proposition claims for the case
$\tmop{Re}z_1>\tmop{Re}z_2>\dots>\tmop{Re}z_n$. By meromorphicity,
the proof is complete.
\end{proof}

\begin{corollary}
For  $w=(n,n+1)$ we have:
\begin{equation} \label{Tlong}
T_{w^{-1}}^*\Lambda_{H,\chi}=(-1)\chi(a_\alpha)c_\alpha(\chi)\Lambda_{H,^w
\chi}.
\end{equation}
For $w=w_\alpha w_{\tilde\alpha}$ where $w_\alpha = (i,i+1)$ with
$1\le i<n$ we have:
\begin{equation} \label{Tshort}
T_{w^{-1}}^*\Lambda_{H,\chi}=(-1)\chi(a_\alpha)
\frac{1-q^{-1}\chi(a_{-\alpha})}{1-q^{-1}\chi(a_\alpha)}
c_\alpha(\chi)c_{\tilde{\alpha}}(\chi)\Lambda_{H,^w
\chi}.\end{equation}
\end{corollary}

More compactly, this can be written:

\emph{For every $w\in \Gamma$,}
\begin{equation} \label{tlambda}
T_{w^{-1}}^*\Lambda_{H,\chi}=(-1)^{l(w)}\prod_{\alpha\in\Phi_{GL}^+,w\alpha<0}
c_\alpha(\chi)\prod_{\alpha\in\Phi_{Sp}^+,w\alpha<0}d_\alpha(\chi)
\end{equation}
\emph{where}
\begin{equation}
d_\alpha(\chi)=\left\{\begin{array}{ll} \chi(a_\alpha) & \textrm{
if }
\alpha \textrm{ is a long root}\\
\chi(a_\alpha)
\frac{1-q^{-1}\chi(a_{-\alpha})}{1-q^{-1}\chi(a_\alpha)} &
\textrm{ if } \alpha \textrm{ is a short root.}\end{array}\right. 
\end{equation}

To complete the proof of our theorem, we bring together equations
(\ref{general})
 and (\ref{tlambda}) to get:

\begin{eqnarray}
\Omega(g_\lambda)= Q^{-1} \prod_{\alpha\in \Phi_{GL}^+} c_\alpha(\chi) \cdot \nonumber \\ \sum_{w\in\Gamma} (-1)^{l(w)} \prod_{\alpha\in
\Phi_{Sp}^+,w\alpha<0} \chi(a_\alpha) \prod_{\alpha\in
\Phi_{Sp}^{S+},
w\alpha<0}\frac{1-q^{-1}\chi(a_{-\alpha})}{1-q^{-1}\chi(a_\alpha)}
{^w\chi}^{-1}\delta^{\frac{1}{2}}(g_\lambda)
\end{eqnarray}
which is just (\ref{nolgroup}) normalized differently.

We now use the fact
that if $k_\alpha$ is an expression which is covariant with
$\alpha$ then
\begin{equation}
\frac{w^{-1}\left(\prod_{\alpha>0}{k_\alpha}\right)}{\prod_{\alpha>0}{k_\alpha}}
=
\frac{\prod_{\alpha>0,w\alpha<0}k_{-\alpha}}{\prod_{\alpha>0,w\alpha<0}
k_\alpha}
\end{equation}
to write the $\chi(a_\alpha)$ factors as:
\begin{equation}
e^{\rho-w^{-1}\rho}(g_\chi)
\end{equation}
and the product over short roots:
\begin{equation}
\frac{\prod_{\alpha\in\Phi_{Sp}^{S+}}(1-q^{-1}e^{w^{-1}\alpha})}
{\prod_{\alpha\in\Phi_{Sp}^{S+}}(1-q^{-1}e^{\alpha})}(g_\chi).
\end{equation}

We also replace $w^{-1}$ by $w w_l$ and $c_\alpha$ by its definition (\ref{cadefinition}), and finally we get:

\begin{eqnarray}
\Omega(g_\lambda)= \frac{\prod_{\alpha\in \Phi_{Sp}^+} (1-q^{-1}e^\alpha)}{Q e^{-\rho}\prod_{\alpha\in\Phi_{GL}^+}(1-e^\alpha)}(g_\chi) \cdot 
\delta^{\frac{1}{2}}(g_\lambda) \mathcal{A} \left(
e^{\rho+\lambda} \prod_{\alpha\in \Phi_{Sp}^{S+}}
(1-q^{-1}e^{-\alpha}) \right)(g_\chi). \nonumber \\ 
\end{eqnarray}

The product \[Q e^{-\rho}\prod_{\alpha\in\Phi_{GL}^+}(1-e^\alpha)\] in the denominator can be ignored, as we could have applied Bernstein's theorem by imposing the normalization:
\begin{equation}\label{finalnorm}
\Lambda(\phi_B)=Q e^{-\rho}\prod_{\alpha\in\Phi_{GL}^+}(1-e^\alpha)(g_\chi)
\end{equation}
instead of (\ref{normalization}). Thus we get the formula of Theorem \ref{mainthm}. Notice that the remaining factor in front of the alternator vanishes exactly when $I(\chi)$ is reducible with the spherical vector generating a proper $G$-subspace. The rest of the expression vanishes identically exactly when $g_\chi$ is non-regular as an element of $\tmop{Sp}_{2n}(\mathbb C)$. In all other cases, the spherical vector generates $I(\chi)$ and is nonzero, hence $\chi$ lies neither on a  singular hypersurface nor on the zero set for the Shalika functional whose existence was provided by Bernstein's theorem, which proves \emph{a posteriori} the existence of such a  non-zero functional. This completes the proof of Theorem \ref{thm}.


\section{A general remark}\label{secgeneral}

A large part of what we did here is virtually independent of the particular setting of the Shalika model, and can be directly transfered to study other unique models induced from closed algebraic subgroups. By slight modifications, one might also hope to include non-algebraic subgroups that are ``twists'' of algebraic ones, but we do not consider those cases here.

Let $G$ be a split reductive algebraic group over $F$ with a fixed integral model over $\mathfrak o$, and let $H$ be any $F$-rational subgroup of $G$ such that the Borel subgroup has a rational open orbit in $G/H$. 
Such subgroups are called spherical and it is known that the quotient $G/H$ has only a finite number of $P$-orbits (cf. \cite{Br} or \cite{V}). For simplicity we will choose a Borel $P$ such that this orbit is represented by the element 1. We set $K=G(\mathfrak o)$, a maximal compact subgroup, and assume that all double $P\backslash G/H$ cosets have a representative in $K$. We also let $B$ be the standard Iwahori, namely the inverse image (in $K$) of $P(\mathfrak{o/p})$ under the reduction map. We also fix a maximal split torus $T\subset P$. Let $\Psi$ be a character (possibly trivial) of $H$, and let us assume that $\Psi$ is trivial on $H_0=H\cup K$. 

The above data define a representation $\tmop{Ind}_H^G(\Psi)$, which has good chances of being multiplicity-free. For the unramified spectrum, one can examine this question by a double $P\backslash G/H$ decomposition. For our purposes, we will just imitate the setting of our present work and \emph{assume} the following:
\begin{enumerate}
\item Let $D$ denote the subvariety of $X(T)$ (the complex torus of unramified characters of the maximal torus of $G$), defined by the condition $\left.\chi\delta_P^{\frac{1}{2}}\right|_{P\cap H}=\left.\Psi \frac{\delta_H}{\delta_{P\cap H}}\right|_{P\cap H}$, as in section \ref{secvanishing}. Assume that the natural integral on the open orbit, representing an intertwining functional $I(\chi)\to \tmop{Ind}_H^G(\Psi)$, namely:
\begin{equation*}
\Lambda(\phi)=\int_{P\cap H\backslash H} \phi(h)\Psi^{-1}(h) dh
\end{equation*}
converges for $\chi$ in some open subset of $D$. (Notice that we have to take all modular characters into account, in contrast to Proposition \ref{propvanishing} where $H$ and $P\cap H$ turned out to be unimodular. If the modular characters of $H$ and $P\cap H$ do not agree on $P\cap H$, this is not a well-defined integral, in which case one has to ``lift'' the distribution to the group and use an expression analogous to (\ref{Delta}.))
\item For generic $\chi$ in the aforementioned open set, the space of intertwining operators: $I(\chi)\to\tmop{Ind}_H^G(\Psi)$ is one-dimensional.
\end{enumerate} 

Under these assumptions, one can use the methods employed here to compute the values of the spherical vector $\Omega_{\chi}\in\tmop{Ind}_H^G(\Psi)$ in the image of $I(\chi)$ \emph{on sufficiently large anti-dominant elements of the torus}, i.e. $g=\lambda^{-1}(\varpi)$, with $\lambda:\mathbb G_m\to T$ a dominant cocharacter such that $|\alpha(g)|$ is sufficiently large for every positive root $\alpha$. In certain cases, like in our example, the method will work for a full class of representatives of $H\backslash G/K$ cosets, but in general it only gives asymptotic results.

More precisely:
Bernstein's theorem guarantees existence, uniqueness (generically) and rationality of the intertwining operators, the general Casselman-Shalika-Hironaka formula (\ref{goal}) still holds, and the Casselman coefficients $a_w$ will vanish, except for $w\in\Gamma$ where $\Gamma$ is the stabilizer (in $W$) of the relation $\left.\chi\delta_P^{\frac{1}{2}}\right|_{P\cap H}=\left.\Psi \frac{\delta_H}{\delta_{P\cap H}}\right|_{P\cap H}$. Moreover, for $g=g_{-\lambda}$, with $\lambda\ge\lambda_0$, $\lambda_0$ a large enough dominant co-character, we will still have $gN_1^-g^{-1}\subset P_0 H_0$, making the computation of proposition \ref{depproposition} valid. The corresponding corollary is, in general, that 
\begin{equation}
\frac{a_w(g_{-\lambda})}{a_w(g_{-\lambda_0})}=\frac{{^w\chi}^{-1} \delta^\frac{1}{2}(g_\lambda)}{{^w\chi}^{-1} \delta^\frac{1}{2}(g_{\lambda_0})}
\end{equation}
Finally, one computes the ``functional equations'' as we did in section \ref{secfe}, but instead of the functions $\phi_{BwB}$ one uses the functions $R_{g_{-\lambda_0}}\phi_{BwB}$. The latter computation will basically involve computing an explicit factorization for $wN$ in the form $P\cdot H$, for $w$ in a set of generators of $\Gamma$.

Hence the problem of computing such an explicit formula - which typically involved a maze of difficult and case-specific considerations - has been reduced to the straightforward computation of this factorization, at least in order to get asymptotic results.


\begin{thebibliography}{}

\bibitem{AG} Ash, Avner; Ginzburg, David. \emph{$p$-adic $L$-functions for ${\rm GL}(2n)$.}  Invent. Math.  116  (1994),  no. 1-3, 27--73. 

\bibitem{BB} Beineke, Jennifer; Bump, Daniel. \emph{A summation formula for divisor functions associated to lattices.} Preprint, 2004.

\bibitem{Be}
Bernstein, Joseph. \emph{Letter to Piatetski-Shapiro.} Fall 1985, unpublished.

\bibitem{Br}  Brion, Michel, \emph{Quelques propri\'et\'es des espaces homog\`enes sph\'eriques.} Manuscripta Math.  55  (1986),  no. 2, 191--198.


\bibitem{BF}  Bump, Daniel; Friedberg, Solomon. \emph{The exterior square automorphic $L$-functions on ${\rm GL}(n)$.}  Festschrift in honor of I. I. Piatetski-Shapiro on the occasion of his sixtieth birthday, Part II (Ramat Aviv, 1989),  47--65, Israel Math. Conf. Proc., 3, Weizmann, Jerusalem, 1990.


\bibitem{BFG} Bump, Daniel; Friedberg, Solomon; Ginzburg, David. \emph{Whittaker-orthogonal models, functoriality, and the Rankin-Selberg method.}  Invent. Math.  109  (1992),  no. 1, 55--96.




\bibitem{Cas}
Casselman, W. \emph{Introduction to the theory of admissible
representations of $p$-adic reductive groups.} Draft, 1 May 1995,
currently available at:
http://www.math.ubc.ca/$\sim$cass/research/p-adic-book.dvi.


\bibitem{C}
Casselman, W. \emph{The unramified principal series of $p$-adic groups. I. The spherical function.}  Compositio Math.  40  (1980), no. 3, 387--406.

\bibitem{CS}
Casselman, W.; Shalika, J. \emph{The unramified principal series of $p$-adic groups. II. The Whittaker function.}  Compositio Math.  41  (1980), no. 2, 207--231.

\bibitem{FJ} Friedberg, Solomon; Jacquet, Herv\'e. \emph{Linear periods.}  J. Reine Angew. Math.  443  (1993), 91--139.


\bibitem{GPS}
Gelbart, S. ; Piatetski-Shapiro, I. \emph{L-functions for $G\times
\tmop{GL}(n)$}, in:  Gelbart, Stephen; Piatetski-Shapiro, Ilya; Rallis, Stephen. \emph{Explicit constructions of automorphic $L$-functions.} Lecture Notes in Mathematics, 1254. Springer-Verlag, Berlin, 1987.

\bibitem{GRS}
Ginzburg, David; Rallis, Stephen; Soudry, David. \emph{Generic automorphic forms on ${\rm SO}(2n+1)$: functorial lift to ${\rm GL}(2n)$, endoscopy, and base change.}  Internat. Math. Res. Notices  2001,  no. 14, 729--764.

\bibitem{H} 
 Hironaka, Yumiko. \emph{Spherical functions and local densities on Hermitian forms.}  J. Math. Soc. Japan  51  (1999),  no. 3, 553--581.


\bibitem{JR}
 Jacquet, Herv\'e; Rallis, Stephen. \emph{Uniqueness of linear periods.}  Compositio Math.  102  (1996),  no. 1, 65--123.

\bibitem{JS}
 Jacquet, Herv\'e; Shalika, Joseph. \emph{Exterior square $L$-functions.}  Automorphic forms, Shimura varieties, and $L$-functions, Vol. II (Ann Arbor, MI, 1988),  143--226, Perspect. Math., 11, Academic Press, Boston, MA, 1990.

\bibitem{MR}
Mao, Zhengyu;  Rallis, Stephen. \emph{Preprint.}

\bibitem{O}
 Offen, Omer. \emph{Relative spherical functions on $p$-adic symmetric spaces (three cases).}  Pacific J. Math.  215  (2004),  no. 1, 97--149. 


\bibitem{Sato}
Sato, Fumihiro. \emph{Enumeration of subgroups of finite abelian $p$-groups,
local densities of square matrices and Fourier coefficients of
Eisenstein series of $\tmop{GL}_n$.} Preprint, 2004.

\bibitem{V}
 Vinberg, \`E. B. \emph{Complexity of actions of reductive groups.} (Russian)  Funktsional. Anal. i Prilozhen.  20  (1986),  no. 1, 1--13, 96.

\end{thebibliography}
\end{document}